\documentclass[12pt]{article}
\usepackage{indentfirst} 
\usepackage{amsfonts,amssymb,amsmath}
\usepackage{bm}          
\usepackage{cite}
\usepackage{url}

\DeclareMathOperator{\sgn}{\rm sgn}

\numberwithin{equation}{section}  
\begin{document}

\bibliographystyle{plain}

\def\today{October 10, 2013 \\[2mm]
           revised February 21, 2014}

\title{Bivariate Generating Functions for a Class of Linear
       Recurrences: General Structure}

\author{
  {\small J. Fernando Barbero G.${}^{1,3}$, Jes\'us Salas${}^{2,3}$, and
          Eduardo J.S. Villase\~nor${}^{2,3}$} \\[4mm]
  {\small\it ${}^1$Instituto de Estructura de la Materia, CSIC} \\[-0.2cm]
  {\small\it Serrano 123, 28006 Madrid, Spain} \\[-0.2cm]
  {\small\tt FBARBERO@IEM.CFMAC.CSIC.ES}               \\[1mm]
  {\small\it ${}^2$Grupo de Modelizaci\'on, Simulaci\'on Num\'erica
                   y Matem\'atica Industrial}  \\[-0.2cm]
  {\small\it Universidad Carlos III de Madrid} \\[-0.2cm]
  {\small\it Avda.\  de la Universidad, 30}    \\[-0.2cm]
  {\small\it 28911 Legan\'es, Spain}           \\[-0.2cm]
  {\small\tt JSALAS@MATH.UC3M.ES, EJSANCHE@MATH.UC3M.ES}  \\[1mm]
  {\small\it ${}^3$Grupo de Teor\'{\i}as de Campos y F\'{\i}sica
             Estad\'{\i}stica}\\[-2mm]
  {\small\it Instituto Gregorio Mill\'an, Universidad Carlos III de
             Madrid}\\[-2mm]
  {\small\it Unidad Asociada al Instituto de Estructura de la Materia, CSIC}
             \\[-2mm]
  {\small\it Madrid, Spain}           \\[-2mm]
  {\protect\makebox[5in]{\quad}}  
  \\
}

\maketitle
\thispagestyle{empty}   

\begin{abstract}
We consider Problem~6.94 posed in the book {\em Concrete Mathematics}
by Graham, Knuth, and Patashnik, and solve it by using bivariate
exponential generating functions. The family of recurrence relations
considered in the problem contains many cases of combinatorial interest for
particular choices of the six parameters that define it. We give a complete
classification of the partial differential equations satisfied by the
exponential generating functions, and solve them in all cases. We also
show that the recurrence relations defining the combinatorial numbers
appearing in this problem display an interesting degeneracy that we study
in detail. Finally, we obtain for all cases the corresponding
univariate row generating polynomials.
\end{abstract}


\clearpage

%
%
%
%
\newcommand{\be}{\begin{equation}}
\newcommand{\ee}{\end{equation}}
\newcommand{\<}{\langle}
\renewcommand{\>}{\rangle}
\newcommand{\widebar}{\overline}
\def\spose#1{\hbox to 0pt{#1\hss}}
\def\ltapprox{\mathrel{\spose{\lower 3pt\hbox{$\mathchar"218$}}
 \raise 2.0pt\hbox{$\mathchar"13C$}}}
\def\gtapprox{\mathrel{\spose{\lower 3pt\hbox{$\mathchar"218$}}
 \raise 2.0pt\hbox{$\mathchar"13E$}}}
\def\textprime{${}^\prime$}
\def\proof{\par\medskip\noindent{\sc Proof.\ }}
\def\qed{\hbox{\hskip 6pt\vrule width6pt height7pt depth1pt \hskip1pt}\bigskip}
\def\proofof#1{\bigskip\noindent{\sc Proof of #1.\ }}
\def\half{\frac{1}{2}}
\def\third{\frac{1}{3}}
\def\twothird{\frac{2}{3}}
\def\smfrac#1#2{\textstyle \frac{#1}{#2}}
\def\smhalf{\smfrac{1}{2} }

\newcommand{\restrict}{\upharpoonright}
\newcommand{\drop}{\setminus}
\renewcommand{\emptyset}{\varnothing}

%
%
\newcommand{\C}{{\mathbb C}}
\newcommand{\D}{{\mathbb D}}
\newcommand{\Z}{{\mathbb Z}}
\newcommand{\N}{{\mathbb N}}
\newcommand{\R}{{\mathbb R}}
\newcommand{\Q}{{\mathbb Q}}

%
%
\newcommand{\TT}{{\mathsf T}}
\newcommand{\HH}{{\mathsf H}}
\newcommand{\VV}{{\mathsf V}}
\newcommand{\JJ}{{\mathsf J}}
\newcommand{\PP}{{\mathsf P}}
\newcommand{\DD}{{\mathsf D}}
\newcommand{\QQ}{{\mathsf Q}}
\newcommand{\RR}{{\mathsf R}}

%
%
\newcommand{\bsigma}{{\boldsymbol{\sigma}}}
\newcommand{\vecbsigma}{{\vec{\boldsymbol{\sigma}}}}
\newcommand{\bpi}{{\boldsymbol{\pi}}}
\newcommand{\vecbpi}{{\vec{\boldsymbol{\pi}}}}
\newcommand{\btau}{{\boldsymbol{\tau}}}
\newcommand{\bphi}{{\boldsymbol{\phi}}}
\newcommand{\bvarphi}{{\boldsymbol{\varphi}}}
\newcommand{\bGamma}{{\boldsymbol{\Gamma}}}

%
%
\newcommand{\psibar}{ {\bar{\psi}} }
\newcommand{\varphibar}{ {\bar{\varphi}} }

%
%
\newcommand{\ba}{ {\bf a} }
\newcommand{\bb}{ {\bf b} }
\newcommand{\bc}{ {\bf c} }
\newcommand{\bp}{ {\bf p} }
\newcommand{\br}{ {\bf r} }
\newcommand{\bs}{ {\bf s} }
\newcommand{\bt}{ {\bf t} }
\newcommand{\bu}{ {\bf u} }
\newcommand{\bv}{ {\bf v} }
\newcommand{\bw}{ {\bf w} }
\newcommand{\bx}{ {\bf x} }
\newcommand{\by}{ {\bf y} }
\newcommand{\bz}{ {\bf z} }
\newcommand{\bone}{ {\mathbf 1} }

%
%
\newcommand{\scra}{{\mathcal{A}}}
\newcommand{\scrb}{{\mathcal{B}}}
\newcommand{\scrc}{{\mathcal{C}}}
\newcommand{\scrd}{{\mathcal{D}}}
\newcommand{\scre}{{\mathcal{E}}}
\newcommand{\scrf}{{\mathcal{F}}}
\newcommand{\scrg}{{\mathcal{G}}}
\newcommand{\scrh}{{\mathcal{H}}}
\newcommand{\scri}{{\mathcal{I}}}
\newcommand{\scrj}{{\mathcal{J}}}
\newcommand{\scrk}{{\mathcal{K}}}
\newcommand{\scrl}{{\mathcal{L}}}
\newcommand{\scrm}{{\mathcal{M}}}
\newcommand{\scrn}{{\mathcal{N}}}
\newcommand{\scro}{{\mathcal{O}}}
\newcommand{\scrp}{{\mathcal{P}}}
\newcommand{\scrq}{{\mathcal{Q}}}
\newcommand{\scrr}{{\mathcal{R}}}
\newcommand{\scrs}{{\mathcal{S}}}
\newcommand{\scrt}{{\mathcal{T}}}
\newcommand{\scru}{{\mathcal{U}}}
\newcommand{\scrv}{{\mathcal{V}}}
\newcommand{\scrw}{{\mathcal{W}}}
\newcommand{\scrx}{{\mathcal{X}}}
\newcommand{\scry}{{\mathcal{Y}}}
\newcommand{\scrz}{{\mathcal{Z}}}

%
%
\newtheorem{theorem}{Theorem}[section]
\newtheorem{definition}[theorem]{Definition}
\newtheorem{proposition}[theorem]{Proposition}
\newtheorem{lemma}[theorem]{Lemma}
\newtheorem{corollary}[theorem]{Corollary}
\newtheorem{conjecture}[theorem]{Conjecture}
\newtheorem{result}[theorem]{Result}
\newtheorem{question}[theorem]{Question}

%
%
\newcommand{\stirlingsubset}[2]{\genfrac{\{}{\}}{0pt}{}{#1}{#2}}
\newcommand{\stirlingcycle}[2]{\genfrac{[}{]}{0pt}{}{#1}{#2}}
\newcommand{\associatedstirlingsubset}[2]{\left\{\!\!%
            \stirlingsubset{#1}{#2} \!\! \right\}}
\newcommand{\associatedstirlingsubsetBis}[2]{\big\{\!\!%
            \stirlingsubset{#1}{#2} \!\! \big\}}
\newcommand{\assocstirlingsubset}[3]{{\genfrac{\{}{\}}{0pt}{}{#1}{#2}}_{\!%
            \ge #3}}
\newcommand{\assocstirlingcycle}[3]{{\genfrac{[}{]}{0pt}{}{#1}{#2}}_{\ge #3}}
\newcommand{\associatedstirlingcycle}[2]{\left[\!\!%
            \stirlingcycle{#1}{#2} \!\! \right]}
\newcommand{\associatedstirlingcycleBis}[2]{\big[\!\!%
            \stirlingcycle{#1}{#2} \!\! \big]}
\newcommand{\euler}[2]{\genfrac{\langle}{\rangle}{0pt}{}{#1}{#2}}
\newcommand{\eulergen}[3]{{\genfrac{\langle}{\rangle}{0pt}{}{#1}{#2}}_{\! #3}}
\newcommand{\eulersecond}[2]{\left\langle\!\! \euler{#1}{#2} \!\!\right\rangle}
\newcommand{\eulersecondBis}[2]{\big\langle\!\! \euler{#1}{#2} \!\!\big\rangle}
\newcommand{\eulersecondgen}[3]%
      {{\left\langle\!\! \euler{#1}{#2} \!\!\right\rangle}_{\! #3}}
\newcommand{\binomvert}[2]{\genfrac{\vert}{\vert}{0pt}{}{#1}{#2}}
\newcommand{\nueuler}[3]{{\genfrac{\langle}{\rangle}{0pt}{}{#1}{#2}}^{\! #3}}

%
%
\section{Introduction} \label{sec.intro}

Graham, Knuth and Patashnik (GKP), in their book {\em Concrete Mathematics}\/
\cite{Graham_94}, posed the following ``research problem''
\cite[Problem~6.94, pp.~319 and~564]{Graham_94}:

\begin{question} \label{problem.GKP}
Develop a general theory of the solutions to the two-parameter recurrence
\begin{equation}
 \binomvert{n}{k} \;=\;
  (\alpha n + \beta k + \gamma)    \, \binomvert{n-1}{k}
  +
  (\alpha' n + \beta' k + \gamma') \, \binomvert{n-1}{k-1} \,+\,
 \delta_{n0}\delta_{k0}
  \label{eq_binomvert}
\end{equation}
for $n,k\in \Z$, assuming that $\binomvert{n}{k} =0$ when $n<0$ or $k<0$.
(Here and in the following $\delta_{ab}$ denotes the Kronecker delta.)
\end{question}

\medskip

Many of the solutions to the recurrence \eqref{eq_binomvert} have been
thoroughly studied in the literature
\cite{Riordan_58,Comtet_74,Graham_94,Bona_07,Bona_12}. They include
classic examples such as the binomial coefficients,
Stirling numbers of several kinds, Eulerian numbers and many others
(see Table~\ref{table.example}).
Particular choices of the parameters defining the problem have been
considered by Neuwirth in \cite{Neuwirth_01}, where he found the solution
of the recursion \eqref{eq_binomvert} for the particular case $\alpha'=0$
by using Galton arrays. Also Spivey \cite{Spivey_11} has found explicit
solutions (using finite differences) for the following three cases:
(S1) $\alpha=-\beta$; (S2) $\beta=\beta'=0$; and
(S3) $\alpha/\beta = \alpha'/\beta' +1$.

The previous studies focused on finding closed expressions for
$\binomvert{n}{k}$ in terms of simpler combinatorial numbers but did not
make significant use of generating functions. After completing the
main computations of this paper, we learned\footnote{
  We thank David Callan for calling our attention to Wilf's paper
  \cite{Wilf_04}, which in turn refers to some earlier work
  by Th\'eor\^et \cite{Theoret_thesis,Theoret_95_1, Theoret_95_2}.
}
that those have been considered in the context of problem \eqref{eq_binomvert}
by Th\'eor\^et \cite{Theoret_thesis,Theoret_95_1, Theoret_95_2} and
Wilf \cite{Wilf_04}. In particular,
Th\'eor\^et finds the exponential generating functions (EGF's) for the
four particular cases explained above, and
Wilf gives a general solution to the partial differential equations (PDE's)
satisfied by the EGF's in terms of hypergeometric functions.
However, in his own words \cite{Wilf_04}:  ``\ldots we obtain
a complete solution also, though its form is very unwieldy''.

In this paper, we study in a systematic way the PDE's satisfied by the
EGF's (in the next formula $P_n(x)$ are the so called
\emph{row generating polynomials})
\begin{equation}
   F(x,y) \;=\; \sum_{n,k\ge 0} \binomvert{n}{k} \, x^k \,
                \frac{y^n}{n!} \;=\; \sum_{n\ge 0}\frac{y^n}{n!}\, P_n(x)\;,
  \label{def_F}
\end{equation}
defined by the sequences of numbers given by the recurrences
\eqref{eq_binomvert}. We propose a classification scheme that leads to a
clean understanding of their solutions.

It is straightforward to show that the EGF \eqref{def_F} associated with the
numbers $\binomvert{n}{k}$ satisfying the recurrence \eqref{eq_binomvert}
is a solution to the PDE
\be
-(\beta + \beta'\, x) \, x  \, F_1 +
(1-\alpha \, y - \alpha' \, x \, y  ) \, F_2
\;=\; (\alpha + \gamma + (\alpha' + \beta' + \gamma')\, x ) \, F \,,
\label{eq_PDE_final}
\ee
with the initial condition $F(x,0)=1$. (Here and in the following $F_i$ denotes the partial derivative of $F$ with respect to its $i$-th variable.)

We classify now the PDE's satisfied by the EGF's solving
Question~\ref{problem.GKP} in terms of
$(\alpha,\beta,\gamma;\alpha',\beta',\gamma')$. The dependence of the
equations on the parameters $(\gamma;\gamma')$ is always fairly simple,
so we will introduce \emph{families} of equations characterized by the
parameters $(\alpha,\beta;\alpha',\beta')$. In the paper we will sometimes
refer to the full set of parameters defining a recurrence
$(\alpha, \beta, \gamma;\alpha',\beta',\gamma')$, and sometimes just to
the family $(\alpha,\beta;\alpha',\beta')$. A careful look at
\eqref{eq_PDE_final} reveals that the two most important parameters are
$\beta$ and $\beta'$. In fact, what really matters is whether these
parameters are zero or non-zero. This leads us to introduce the following
four different \emph{types} of equations:

\begin{definition} \label{def_types}
The PDE's for the EGF's relevant to solve Question~\ref{problem.GKP} are
classified in four different types: {\bf Type~I}: $\beta\beta' \neq 0$;
{\bf Type~II}: $\beta\neq 0$ and $\beta' = 0$;
{\bf Type~III}: $\beta =0$ and $\beta' \neq 0$; and
{\bf Type~IV}:  $\beta = \beta' = 0$.
\end{definition}

\noindent \textbf{Remarks:}
1. It is important to notice that although  the parameters
$(\alpha,\beta,\gamma;\alpha',\beta',\gamma')$ uniquely determine the numbers  
$\binomvert{n}{k}$, the converse is \emph{not} true. For example,
the \emph{trivial sequence} $\binomvert{n}{k}= \delta_{n0}\delta_{k0}$
can be obtained through equations of any type by choosing the parameters to
be $(\alpha,\beta,-\alpha;\alpha',\beta',-\alpha'-\beta')$, regardless of the
specific values of $\alpha$, $\beta$, $\alpha'$, and $\beta'$, as can be
easily seen by looking at Eq.~\eqref{eq_PDE_final}. We will explore this
phenomenon in the present paper and identify all the possible indeterminacies
of this type. The reason why we classify equations instead of their solutions
is a direct consequence of this fact.

2. If the numbers $\binomvert{n}{k}$ satisfy a recursion of the form
\eqref{eq_binomvert} for certain  parameters
$(\alpha,\beta,\gamma;\alpha',\beta',\gamma')$, then the numbers
\be
\binomvert{n}{k}^\star \;=\;  \binomvert{n}{n-k}
\ee
also satisfy a recursion of the form \eqref{eq_binomvert} with parameters
\be
(\alpha^\star,\beta^\star,\gamma^\star;
 \alpha'^\star,\beta'^\star,\gamma'^\star) \;=\;
(\alpha'+\beta',-\beta',\gamma';\alpha+\beta,-\beta,\gamma)\,.
\label{def_duality}
\ee
This involution was already introduced by Th\'eor\^et
\cite[Eq.~(34)]{Theoret_95_1}.
Hence, every Type~II family is the $\star$-image of a Type~III family
and vice versa. On the other hand, the classes of Types~I
and~IV are both closed under the $\star$-map, and it makes sense to
talk about self-dual (symmetric) numbers as those satisfying
\be
\binomvert{n}{k} \;=\; \binomvert{n}{n-k} \;=\; \binomvert{n}{k}^\star\,.
\ee
The families $(\alpha,\beta,\gamma;\alpha+\beta,-\beta,\gamma)$ are,
for example, self-dual.

A more general involution can be defined as follows: let us denote by
$\bm{z}=(x,y)$ the variables of the EGF \eqref{def_F}, and by
$\bm{\mu}=(\alpha,\beta,\gamma;\alpha',\beta',\gamma')$ the parameters of the
corresponding recurrence \eqref{eq_binomvert}, so the EGF can be compactly
rewritten as $F(\bm{z};\bm{\mu})$. We now define the $\star$-image of $F$ as:
\be
F^\star(\bm{z};\bm{\mu}) \;=\; F\left(\mathcal{M}_1(\bm{z},\bm{\mu});
                                      \mathcal{M}_2(\bm{\mu}) \right) \,,
\ee
where $\mathcal{M}_2(\bm{\mu})$ is an involution (i.e.,
$\mathcal{M}_2(\mathcal{M}_2(\bm{\mu}))=\bm{\mu}$), and the function
$\mathcal{M}_1$ satisfies
\be
\mathcal{M}_1\left( \mathcal{M}_1(\bm{z},\bm{\mu}),
                    \mathcal{M}_2(\bm{\mu}) \right) \;=\; \bm{z}\,,
\qquad \text{for all $\bm{z}$.}
\ee
Then, this $\star$-map is obviously an involution: i.e.,
$(F^{\star})^\star(\bm{z};\bm{\mu})=F(\bm{z};\bm{\mu})$.
The simplest of these involutions
(i.e., those for which $\mathcal{M}_1$ is only a function of $\bm{z}$)
are listed in Table~\ref{table.involutions}. In Ref.~\cite{BSV} we will
discuss in detail a more involved case.

\medskip

%
%
\begin{table}
\small
\begin{tabular}{cllll}
\hline\hline\\[-3mm]
Type & \multicolumn{1}{c}{Family} & $(\gamma,\gamma')$ & Description &
       \multicolumn{1}{c}{Entry} \\[1mm]
\hline \\[-3mm]
I    &  $(0,1;1,-1)$ & $(1,0)$ & Eulerian numbers $\euler{n}{k}$
        \cite{Graham_94} & A173018\\
     &               &         & \phantom{Eulerian numbers $\euler{n}{k}$}
        \cite{Riordan_58,Comtet_74} & A008292\\[1mm]
\cline{2-5}\\[-3mm]
     &  $(0,1;2,-1)$ & $(1,-1)$ & Second-order Eulerian numbers
        $\eulersecond{n}{k}$ \cite{Gessel_78,Graham_94} & A008517\\[1mm]
\cline{2-5}\\[-3mm]
     &  $(0,1;3,-1)$ & $(1,-2)$ & Third-order Eulerian numbers
        \cite{Park_94,Janson_11} & A219512\\[1mm]
\cline{2-5}\\[-3mm]
    &  $(0,1;\nu,-1)$ & $(1,1-\nu)$ & $\nu$-order Eulerian numbers
       $\nueuler{n}{k}{\nu}$ \cite{BSV}&        \\[1mm]
\cline{2-5}\\[-3mm]
     &  $(0,1;0,1)$  & $(0,0)$  & $\mathrm{Surj}(n,k)$ \cite{Fekete_94,Bona_07}
                     & A019538 \\[1mm]
\cline{2-5}\\[-3mm]
     &  $(0,1;1,1)$ & $(0,-1)$ & Ward numbers
    $\associatedstirlingsubsetBis{n+k}{k}$ \cite{Fekete_94} & A134991 \\[1mm]
\cline{2-5}\\[-3mm]
 &  $(0,1;\nu,1)$ & $(0,-\nu)$ & $\nu$-order Ward numbers \cite{BSV} & \\[1mm]
\cline{2-5}\\[-3mm]
     &  $(1,1;1,1)$ & $(-1,-1)$ &
        $\associatedstirlingcycleBis{n+k}{k}$ \cite{Fekete_94} &      \\[1mm]
\hline\\[-3mm]
II   & $(0,1;0,0)$ & $(0,1)$ & Stirling subset numbers
       $\stirlingsubset{n}{k}$ \cite{Graham_94} & A008277 \\[1mm]
\cline{2-5}\\[-3mm]
     &  $(-1,-1;0,0)$ & $(1,-1)$ & Lah numbers
        $L_{n,k}$ \cite{Riordan_58,Comtet_74} & A008297\\[1mm]
\cline{2-5}\\[-3mm]
     &  $(1,1;0,0)$ & $(-1,1)$ & Unsigned Lah numbers
        $L(n,k)$ \cite{Wagner_96} & A105278\\[1mm]
\cline{2-5}\\[-3mm]
     &  $(2,1;0,0)$ & $(-2,1)$ & Generalization of $\stirlingsubset{n}{k}$
        and $L(n,k)$ \cite{Lang_00,Lang_09} & A035342\\[1mm]
\cline{2-5}\\[-3mm]
     & $(1,-1;0,0)$ & $(0,1)$ & $[n\ge k]\, n!/k!$ \cite{Graham_94,Walsh_11} &
       A094587 \\[1mm]
\cline{2-5}\\[-3mm]
     & $(r-1,1;0,0)$ & $(1-r,1)$ & $S(r;n,k)$ \cite{Lang_00,Lang_09} &\\[1mm]
\hline\\[-3mm]
III  & $(0,0;0,1)$ & $(1,0)$ & $\mathrm{Inj}(n,k)$ \cite{Fekete_94} &
       A008279 \\[1mm]
\cline{2-5}\\[-3mm]
     & $(0,0;-2,1)$ & $(1,0)$ & Coefficients of Laguerre polynomials & \\
     &              &         & in reverse order \cite{Arfken_05}
     & A021010  \\[1mm]
\cline{2-5}\\[-3mm]
     & $(0,0;-1,1)$ & $(1,0)$ & $(-1)^k \stirlingsubset{n}{n-k}$ & A106800
     \\[1mm]
\cline{2-5}\\[-3mm]
     & $(1,0;1,1)$ & $(-1,-1)$ & Ramanujan function $Q_{n+1,k}(-1)$
       \cite{Dumont_96,Zeng_99} & A075856 \\[1mm]
\cline{3-5}\\[-3mm]
     &             & $(0,-2)$ & Ramanujan function $Q_{n,k}(1)$ \cite{Zeng_99}
     & A217922 \\[1mm]
\cline{3-5}\\[-3mm]
     &             & $(0,-1)$ & Ramanujan function $Q_{n+1,k}(0)$
    \cite{Zeng_99} & A054589 \\[1mm]
\hline\\[-3mm]
IV   & $(0,0;0,0)$ & $(1,1)$ & Binomial coefficients $\binom{n}{k}$
        \cite{Graham_94} & A007318 \\[1mm]
\cline{2-5}\\[-3mm]
     & $(1,0;0,0)$ & $(-1,1)$ & Stirling cycle numbers $\stirlingcycle{n}{k}$
        \cite{Graham_94} & A132393 \\[1mm]
\cline{2-5}\\[-3mm]
     & $(-1,0;0,0)$ & $(1,1)$ & Stirling numbers of the 1st kind $s(n,k)$
        \cite{Comtet_74,Bona_07,Bona_12} & A008275 \\[1mm]
\hline\hline
\end{tabular}
\caption{\label{table.example}
Some sequences of combinatorial interest satisfying \eqref{eq_binomvert}.
For each sequence, we give the type of the PDE for its EGF, the parameters
defining its family $(\alpha,\beta;\alpha',\beta')$, the coefficients
$(\gamma,\gamma')$, its description, and the corresponding entry in
Ref.~\cite{Sloane} (if any). More examples can be found in similar
tables in Refs.~\cite{Theoret_thesis,Theoret_95_1}.
}
\end{table}

%
%
\begin{table}
\small
\begin{tabular}{lll}
\hline\hline\\[-3mm]
Involution & Parameter transformation & EGF transformation\\[1mm]
\hline\\[-3mm]
$\displaystyle\binomvert{n}{k}\rightarrow\binomvert{n}{n-k}$&
\begin{minipage}{6.5cm} $(\alpha,\beta,\gamma;\alpha',\beta',\gamma')
\rightarrow\\\hspace*{1cm}
(\alpha'+\beta',-\beta',\gamma';\alpha\!+\beta,-\beta,\gamma)$\end{minipage}
\vspace*{1mm}&$\displaystyle F^*(x,y)=F\left(1/x,xy\right)$\\[1mm]
\hline\\[-3mm]
$\displaystyle\binomvert{n}{k}\rightarrow(-1)^k\binomvert{n}{n-k}$&
\begin{minipage}{6.5cm} $(\alpha,\beta,\gamma;\alpha',\beta',\gamma')
\rightarrow\\\hspace*{1cm}
(\alpha'+\beta',-\beta',\gamma';-\alpha\!-\beta,\beta,-\gamma)$\end{minipage}
\vspace*{1mm}&
$\displaystyle F^*(x,y)=F\left(-1/x,-xy\right)$\\[1mm]
\hline\\[-3mm]
$\displaystyle\binomvert{n}{k}\rightarrow(-1)^k\binomvert{n}{k}$&
\begin{minipage}{6.5cm} $(\alpha,\beta,\gamma;\alpha',\beta',\gamma')
\rightarrow\\\hspace*{1cm}(\alpha,\beta,\gamma;-\alpha',-\beta',-\gamma')$
\end{minipage}\vspace*{1mm}&$F^*(x,y)=F(-x,y)$\\[1mm]
\hline\\[-3mm]
$\displaystyle\binomvert{n}{k}\rightarrow(-1)^{n-k}\binomvert{n}{k}$&
\begin{minipage}{6.5cm} $(\alpha,\beta,\gamma;\alpha',\beta',\gamma')
\rightarrow\\\hspace*{1cm}(-\alpha,-\beta,-\gamma;\alpha',\beta',\gamma')$
\end{minipage}\vspace*{1mm}&$F^*(x,y)=F(-x,-y)$\\[1mm]\hline\\[-3mm]
$\displaystyle\binomvert{n}{k}\rightarrow(-1)^n\binomvert{n}{k}$&
\begin{minipage}{6.5cm} $(\alpha,\beta,\gamma;\alpha',\beta',\gamma')
\rightarrow\\\hspace*{1cm}(-\alpha,-\beta,-\gamma;-\alpha',-\beta',-\gamma')$
\end{minipage}\vspace*{1mm}&$F^*(x,y)=F(x,-y)$\\[1mm]
\hline\hline
\end{tabular}
\caption{\label{table.involutions}
Involutions, their effects on the parameters
$(\alpha,\beta,\gamma;\alpha',\beta',\gamma')$, and on the EGF's.
They can be immediately checked by performing the change of variables
defining the EGF transformation in the PDE \eqref{eq_PDE_final}.}
\end{table}

The plan of the paper is the following. After this introduction,
Section~\ref{sec.solutionF} is devoted to the derivation of
the EGF's for the four distinct types of equations (cf.
Definition~\ref{def_types}). In Section~\ref{sec.ambiguities}, by using
generating functions, we classify the parameter ambiguities in the problem;
i.e., the possibility of obtaining the same solution to the recurrence
with different sets of parameters. Finally, Section~\ref{sec.poly} is
devoted to the study of polynomial generating functions in one variable,
and Appendix~\ref{sec.particular.cases} compiles some particular cases
for which their EGF's can be computed in closed form.

%
%
\section{Exponential generating functions}
\label{sec.solutionF}

As we mentioned in the introduction, Wilf \cite{Wilf_04} has given general
solutions for the EGF's solving \eqref{eq_PDE_final} in terms of
hypergeometric functions. This type of solution partially hides the
structure of the EGF's, and it does not make it easy to see which parameter
choices are the most relevant ones. To avoid these problems we have
introduced a classification of the PDE's for the EGF's in four types
that we discuss in the following subsections.

%
%
\subsection{Type I equations}
\label{sec.typeI}

In this case both $\beta$ and $\beta'$ are non-zero. This fact can be used
to rewrite the PDE \eqref{eq_PDE_final} in a simpler way by performing the
change of variables $(x,y) \mapsto (X,Y)$ defined by
\begin{subequations}
\label{def_XY}
\begin{align}
X &\;=\; \left| \frac{\beta'}{\beta} \right| \, x
  \;=\; \sigma \, \frac{\beta'}{\beta} \, x \,, \qquad
  \text{where $\sigma = \sgn(\beta\beta')$,} \label{def_X} \\
Y &\;=\; \beta \, y  \label{def_Y} \,.
\end{align}
\end{subequations}
Notice that the signs of $x$ and $X$ are the same. This is not strictly
necessary, but it is convenient to avoid absolute values in $\log x$ or
$\log X$, as we will consider $x>0$ and, hence, $X>0$.
The function $F(x,y)$ is then given by a function $\mathcal{F}(X,Y)$ via
the relation
\be
F(x,y) \;=\; \mathcal{F}(X,Y) \;=\; \mathcal{F}\left(
      \sigma \, \frac{\beta'}{\beta} \, x, \beta \, y \right) \,,
\label{def_Fcal}
\ee
where $\mathcal{F}(X,Y)$ satisfies the PDE
\be
 - (1 + \sigma \, X)\, X  \, \mathcal{F}_1 + (1- r\, Y - \sigma \, r' \,X \, Y) \, \mathcal{F}_2 \;=\; (s - \sigma \, s'\, X) \, \mathcal{F}\,,
\label{eq_PDE_XY_final}
\ee
with initial condition
\be
\mathcal{F}(X,0) \;=\; 1 \,.
\label{eq_PDE_F_XY_initial}
\ee
The parameters $r$, $r'$, $s$, $s'$ appearing in \eqref{eq_PDE_XY_final}
are defined as:
\be
r \;=\; \frac{\alpha}{\beta}\,, \quad
r'\;=\; \frac{\alpha'}{\beta'}\,, \quad
s \;=\; \frac{\alpha+\gamma}{\beta}\,, \quad
s'\;=\; -1 - \frac{\alpha'+\gamma'}{\beta'} \,.
\label{def_rs}
\ee
We have, hence, reduced the number of continuous parameters in two units:
from $(\alpha,\beta,\gamma;\alpha', \beta', \gamma')$ to  $(r, s;r', s')$,
plus a discrete parameter $\sigma=\pm1$). By doing this
the expression for the PDE \eqref{eq_PDE_XY_final} becomes simpler
than the original one \eqref{eq_PDE_final}. By using now the well-known method
of characteristics, it is straightforward to prove the following
\begin{theorem} \label{theo.typeI}
The solution $\mathcal{F}$ to \eqref{eq_PDE_XY_final} satisfying
$\mathcal{F}(X,0) \;=\; 1$  is given by
\begin{multline}
\mathcal{F}(X,Y) \;=\; \left(
\frac{ G_{r,r',\sigma}\big(Y\, X^{-r} \, (1+\sigma\, X)^{r-r'} +
       G_{r,r',\sigma}^{-1}(X)\big) }{X} \right)^s \\ 
   \times
\left( \frac{ 1 + \sigma\, X}{1 + \sigma\, G_{r,r',\sigma}
       \big(Y\, X^{-r} \, (1+\sigma\, X)^{r-r'} + G_{r,r',\sigma}^{-1}(X)\big)}
       \right)^{s+s'}\,,
\label{FI}
\end{multline}
where
\be
G_{r,r',\sigma}^{-1}(X) \;=\;  \sum\limits_{k\in\Z_0 \setminus\{r\}}
  \sigma^k \, \binom{-1-r'+r}{k} \, \frac{X^{k-r}}{k-r} +
  \chi_{\Z_0}(r) \, \sigma^r \, \binom{-1-r'+r}{r} \, \log X
\label{def_ginverse_typeI}
\ee
for $0< X < 1$, $\chi_A$ denotes the characteristic function of the set $A$,
and $\Z_0 = \N \cup \{0\}$.
\end{theorem}

It is important to point out here that in many cases of combinatorial
interest, the series defining $G_{r,r',\sigma}^{-1}$ can be summed in closed
form to give simple functions (see Appendix~\ref{sec.particular.cases}).
Something similar happens for the other types of equations.

%
%
\subsection{Type II equations}
\label{sec.typeII}

This type corresponds to $\beta\neq0$ and $\beta'=0$; then the PDE
\eqref{eq_PDE_final} simplifies to
\be
 - \beta \, x \, F_1 +(1 - (\alpha + \alpha'\, x)\, y)\, F_2\;=\;
(\alpha + \gamma + (\alpha' + \gamma')\, x ) \, F \,.
\label{eq_PDE_betap_zero}
\ee
This equation can be solved by using the method of characteristics.
The solution is the content of the following theorem.
\begin{theorem} \label{theo.typeII}
When $\beta\neq 0$ and $\beta'=0$,  the EGF is given by
\begin{multline}
F(x,y) \;=\; \left(
\frac{  G_{\alpha,\beta,\alpha'}\left(y \, x^{-\alpha/\beta}
        e^{-\alpha' x/\beta} + G_{\alpha,\beta,\alpha'}^{-1}(x) \right) }{x}
        \right)^{(\alpha+\gamma)/\beta}  \\
\times \exp\left[ -
  \frac{\alpha'+\gamma'}{\beta}\,
  \left(x -  G_{\alpha,\beta,\alpha'}\left(y \, x^{-\alpha/\beta} \,
  e^{-\alpha' x/\beta} + G_{\alpha,\beta,\alpha'}^{-1}(x) \right)\right)
  \right], 
\label{def_gf_betap=0}
\end{multline}
where $G_{\alpha,\beta,\alpha'}^{-1}(x)$ is defined for any $x>0$ as:
\be
G_{\alpha,\beta,\alpha'}^{-1}(x) \;=\;
    \sum\limits_{k\in \Z_0\setminus\{\alpha/\beta\}}
    \frac{ (-\alpha'/\beta)^k}{k! \, \beta}  \,
    \frac{x^{k-\alpha/\beta}}{k-\alpha/\beta} +
    \chi_{\Z_0}(\alpha/\beta) \,
    \frac{(-\alpha'/\beta)^{\alpha/\beta}}{(\alpha/\beta)!\, \beta }
          \, \log x\,.
\label{def_ginverse_typeII}
\ee
\end{theorem}

%
%
\subsection{Type III equations}
\label{sec.typeIII}

This type corresponds to $\beta=0$ and $\beta'\neq0$. The PDE
\eqref{eq_PDE_final} reads now
\be
-\beta' \, x^2 F_1 +(1 - (\alpha + \alpha'\, x)\, y) \, F_2 \;=\;
(\alpha + \gamma + (\alpha' +\beta'+ \gamma')\, x ) \, F
\label{eq_PDE_beta=zero}
\ee
and we have the following
\begin{theorem} \label{theo.typeIII}
When $\beta= 0$ and $\beta'\neq 0$, the EGF is given by
\begin{multline}
F(x,y) \;=\; \left(
  \frac{ G_{\alpha,\alpha',\,\beta'}\left( y \, x^{-\alpha'/\beta'} \,
         e^{\alpha/(\beta' x)} + G_{\alpha,\alpha',\,\beta'}^{-1}(x)
         \right)}{x} \right)^{1+(\alpha'+\gamma')/\beta'} \\
 \times \exp\left[ \frac{\alpha+\gamma}{\beta'}\,
  \left[\frac{1}{x} -
  \frac{1}{G_{\alpha,\alpha',\,\beta'}\left( y \, x^{-\alpha'/\beta'} \,
           e^{\alpha/(\beta' x)} + G_{\alpha,\alpha',\,\beta'}^{-1}(x)
  \right)}\right] \right],
\label{def_gf_beta=0}
\end{multline}
where $G_{\alpha,\alpha',\,\beta'}^{-1}(x)$ is defined for $x > 0$ as:
\begin{multline}
G_{\alpha,\alpha',\beta'}^{-1}(x) \;=\;  
     -\sum\limits_{k\in \Z_0\setminus\{-1-\alpha'/\beta'\}}
      \frac{(\alpha/\beta')^k}{k!\, \beta'}  \,\frac{1}{k+1+\alpha'/\beta'}\,
      \frac{1}{x^{k+1+\alpha'/\beta'}}  \\
      + \chi_{\N}(-\alpha'/\beta') \,
    \frac{(\alpha/\beta')^{-1-\alpha'/\beta'}}{(-1-\alpha'/\beta')! \, \beta'}
     \, \log x \,.
\label{def_ginverse_typeIII}
\end{multline}
\end{theorem}
It is interesting to note here that the involutions of the first two rows of
Table~\ref{table.involutions} turn equations of Type~II into equations of
Type~III (and vice versa) when the parameters and the arguments of
the EGF's are transformed according to the rules given in any of those rows.
Hence it would have sufficed, in principle, to discuss one of the two types
of equations. We have considered both here for the sake of clarity.

%
%
\subsection{Type IV equations}
\label{sec.typeIV}

This type is characterized by $\beta=\beta'=0$,
and corresponds to case (S2) of Spivey \cite{Spivey_11}.
For the families $(\alpha,0;\alpha',0)$, Eq.~\eqref{eq_PDE_final}
simplifies to the ordinary differential equation
\be
(1 - (\alpha + \alpha'\, x)\, y)\, F_2 \;=\; (\alpha + \gamma + (\alpha' +
     \gamma')\, x ) \, F \,.
\label{eq_PDE_betas=zero}
\ee
Hence, in this case, it is trivial to obtain its closed form solutions
satisfying the initial condition $F(x,0)=1$
(see also \cite[Eq.~(20)]{Theoret_95_1}):

\begin{theorem} \label{theo.typeIV}
When $\beta=\beta'=0$ the EGF is
\be
F(x,y) \;=\; \begin{cases}
\left(1 -(\alpha + \alpha'\, x)\, y\right)^{-\frac{\alpha + \gamma +
   (\alpha'+\gamma')\, x}{\alpha + \alpha' \, x}} &
   \text{if $(\alpha,\alpha')\neq (0,0)$,} \\[2mm]
\exp\left( (\gamma+\gamma'\, x)\, y\right) &
   \text{if $(\alpha,\alpha')= (0,0)$.}
\end{cases}
\label{def_gf_betas=0}
\ee
\end{theorem}
In particular, the EGF's for the Type~IV self-dual families
$(\alpha,0,\gamma;\alpha,0,\gamma)$ are:
\be
F(x,y) \;=\; \begin{cases}
\left(1 -\alpha \,(1 +  x)\, y\right)^{-\frac{\alpha + \gamma}{\alpha} } &
                           \text{if $\alpha\neq 0$,} \\[2mm]
\exp(\gamma\,(1+ x)\, y) & \text{if $\alpha= 0$.}
\end{cases}
\label{def_gf_betas=0_selfdual}
\ee

%
%
\section{Parameter ambiguities} \label{sec.ambiguities}

We discuss here the possibility of obtaining the same combinatorial numbers
$\binomvert{n}{k}$ with different choices of parameters
$(\alpha,\beta, \gamma; \alpha',\beta', \gamma')$. A straightforward way to
do this is to consider Eq.~\eqref{def_F} for the same EGF $F(x,y)$
and \emph{two different sets of parameters}
$(\alpha_1,\beta_1,\gamma_1;\alpha'_1,\beta'_1, \gamma'_1)$ and
$(\alpha_2,\beta_2,\gamma_2;\alpha'_2,\beta'_2, \gamma'_2)$, i.e.,
\begin{subequations}
\label{equations_param_1_2}
\begin{align}
-(\beta_1+\beta'_1 x)x F_1+(1-\alpha_1 y-\alpha'_1 x y)F_2 &\;=\; 
   \left(\alpha_1+\gamma_1+(\alpha'_1+\beta'_1+\gamma'_1)x\right)F\,, 
\label{equation_param_1}\\
-(\beta_2+\beta'_2 x)x F_1+(1-\alpha_2 y-\alpha'_2 x y)F_2 &\;=\;
   \left(\alpha_2+\gamma_2+(\alpha'_2+\beta'_2+\gamma'_2)x\right)F\,,
\label{equation_param_2}
\end{align}
\end{subequations}
with $F(x,0)=1$. By subtracting both equations, we see that a necessary
condition that the EGF $F(x,y)$ must satisfy in order to give rise to the
\emph{same} family of combinatorial numbers is
\begin{equation}
-(\beta_{12}+\beta'_{12} x)x F_1-(\alpha_{12} +\alpha'_{12} x )y F_2
\;=\; \left(\alpha_{12}+\gamma_{12}+(\alpha'_{12}+\beta'_{12}+\gamma'_{12})x
\right)F\,,
\label{eq_fundamental}
\end{equation}
where $\alpha_{12}=\alpha_{1}-\alpha_{2}$, $\beta_{12}=\beta_{1}-\beta_{2}$,
$\gamma_{12}=\gamma_{1}-\gamma_{2}$, $\alpha'_{12}=\alpha'_{1}-\alpha'_{2}$,
$\beta'_{12}=\beta'_{1}-\beta'_{2}$, and $\gamma'_{12}=\gamma'_{1}-\gamma'_{2}$.

The simplest type of ambiguity occurs when $F(x,y)=F(x,0)=1$ (corresponding
to the trivial case for which $\binomvert{n}{k}=\delta_{n0}\delta_{k0}$).
In this case, Eqs.~\eqref{equations_param_1_2}
imply that $\alpha_1+\gamma_1=\alpha_2+\gamma_2=0$ and
$\alpha'_1+\beta'_1+\gamma'_1=\alpha'_2+\beta'_2+\gamma'_2=0$. This means
that any choice of parameters $(\alpha,\beta, \gamma; \alpha',\beta', \gamma')$
such that $\alpha+\gamma=0$ and $\alpha'+\beta'+\gamma'=0$ defines the same
(trivial) family of numbers $\binomvert{n}{k}$. In the following we will
assume that $F(x,y)$ is not constant.

The solutions to Eq.~\eqref{eq_fundamental} are remarkably simple,
and can be written in terms of elementary functions. In order to solve it,
one has to separately consider four cases again:
(i) $\beta_{12}\beta'_{12}\neq0$,
(ii) $\beta_{12}=0$, $\beta'_{12}\neq0$,
(iii) $\beta_{12}\neq0$, $\beta'_{12}=0$, and
(iv) $\beta_{12}=0$, $\beta'_{12}=0$. The solutions for $F(x,y)$ are,
respectively:

\begin{itemize}

\item[(i)] $\displaystyle F(x,y)=
 x^{-\frac{\alpha_{12}+\gamma_{12}}{\beta_{12}}}
 (\beta_{12}+\beta'_{12}x)^{-1+\frac{\alpha_{12}+\gamma_{12}}{\beta_{12}}
 -\frac{\alpha'_{12}+\gamma'_{12}}{\beta'_{12}}}
  \Psi_1\left(yx^{-\frac{\alpha_{12}}{\beta_{12}}}
       (\beta_{12}+\beta'_{12}x)^{\frac{\alpha_{12}}{\beta_{12}}-
        \frac{\alpha'_{12}}{\beta'_{12}}}\right)$.

\item[(ii)] $\displaystyle F(x,y)=\exp\left(
  \frac{\alpha_{12}+\gamma_{12}}{\beta'_{12}}\frac{1}{x}\right)
  x^{-1-\frac{\alpha'_{12}+\gamma'_{12}}{\beta'_{12}}}
  \Psi_2\left(ye^{\frac{\alpha_{12}}{\beta'_{12}}\frac{1}{x}}
              x^{-\frac{\alpha'_{12}}{\beta'_{12}}} \right)$.

\item[(iii)] $\displaystyle F(x,y)=\exp\left(
  -\frac{\alpha'_{12}+\gamma'_{12}}{\beta_{12}}x\right)
  x^{-\frac{\alpha_{12}+\gamma_{12}}{\beta_{12}}}
  \Psi_3\left(ye^{-\frac{\alpha'_{12}}{\beta_{12}}x}
  x^{-\frac{\alpha_{12}}{\beta_{12}}} \right)$.

\item[(iv)] $\displaystyle F(x,y)=y^{-\frac{\alpha_{12}+\gamma_{12}
     +(\alpha'_{12}+\gamma'_{12})x}{\alpha_{12}+\alpha'_{12}x}}\Psi_4(x)$,
     if $\alpha_{12}\neq0$ or $\alpha'_{12}\neq0$.
     If $\alpha_{12}=\alpha'_{12}=0$ then either $\gamma_{12}=\gamma'_{12}=0$
     --and $F(x,y)$ is arbitrary-- or we must have $F(x,y)=0$.

\end{itemize}
The functions $\Psi_j$, $j=1,\ldots,4$, appearing in the preceding
expressions are arbitrary at this stage.

By demanding that $F(x,0)=1$ for all $x>0$ in a neighborhood of $0$,
and requiring that Eq.~\eqref{equation_param_1} (or equivalently,
\eqref{equation_param_2}) is satisfied, we get the conditions that the
parameters must satisfy in order to have non-trivial parameter
indeterminacies, and also the functional forms of the functions $\Psi_j$.
Once the EGF's are obtained, it is easy to derive closed formulas for the
combinatorial numbers that they encode. All these steps are straightforward,
so we just give here the final form of the degenerate families:

\begin{itemize}
\item $\displaystyle (\alpha,\beta, \gamma; \alpha',\beta', \gamma')=
       \left(\alpha,\alpha+\rho G,G-\alpha;
        -\rho H, \rho H+\alpha\frac{H}{G},H-\alpha\frac{H}{G}\right)$. Then
\be
\binomvert{n}{k}_{\alpha,G,H} \;=\;
\binom{n}{k}\, \left( \frac{H}{G}\right)^k \, \prod_{j=0}^{n-1}(G+\alpha j)\,,
\ee
independent of $\rho$. Notice that $G=\alpha+\gamma$ and
$H=\alpha'+\beta'+\gamma'$.

\item $\displaystyle (\alpha,\beta, \gamma; \alpha',\beta', \gamma')=
       (\alpha,-\alpha,-\alpha;-\alpha', L-\alpha',\gamma')$. Then
\be
\binomvert{n}{k}_{L,\gamma'}\;=\; \delta_{nk} \, \prod_{j=1}^n(\gamma'+L j)\,,
\ee
independent of $\alpha,\alpha'$. In this case $L=\alpha'+\beta'$.

\item $\displaystyle (\alpha,\beta, \gamma; \alpha',\beta', \gamma')=
       (\alpha,\beta,M-\alpha;0,\beta',-\beta')$. Then
\be
\binomvert{n}{k}_{M,\alpha} \;=\; \delta_{k0} \,
\prod_{j=0}^{n-1}(M+\alpha j)\,,
\ee
independent of $\beta,\beta'$. In this case $M=\alpha+\gamma$.

\item $\displaystyle (\alpha,\beta, \gamma; \alpha',\beta', \gamma')=
       (\alpha,\beta,-\alpha;\alpha', \beta',-\alpha'-\beta')$.
\be
\binomvert{n}{k} \;=\; \delta_{k0}\delta_{n0}\,,
\ee
independent of $\alpha,\beta,\alpha',\beta'$. This is the trivial case.
\end{itemize}

As we can see, by adjusting the parameters within each of these families,
it is possible to change the type of the PDE for their EGF. This is
the reason why we introduced a classification for the equations instead
of the combinatorial numbers themselves.

%
%
\section{Polynomial generating functions in one variable} \label{sec.poly}

We study here the one-variable polynomials
$P_n(x)$ defined in \eqref{def_F} for the four types of equations used
in the solution to
Question~\ref{problem.GKP}, as defined in Definition~\ref{def_types}.
We get these polynomials from the corresponding EGF's $F(x,y)$ by
employing complex-variable
methods. We will start with the most general case (Type~I), and will work out
the proof of the main theorem with some detail. For the other cases,
the corresponding proofs are very similar, so we will just sketch them
for the sake of brevity.

%
%
\subsection{Type I case} \label{sec.poly_typeI}

When $\beta\beta'\neq 0$, it is convenient to work with the variables
$X,Y$ introduced in \eqref{def_XY}, and define the auxiliary polynomials
\be
\mathcal{P}_n(X) \;=\; n!\, [Y^n]\mathcal{F}(X,Y)\,,
\label{def_Pcal}
\ee
so that
\be
P_n(x) \;= \;\beta^n  \;
\mathcal{P}_n\left(\sigma\frac{\beta^\prime}{\beta}x\right)\,.
\label{def_Pn_vs_Pcal}
\ee

In this section we will get concrete expressions for  $\mathcal{P}_n(X)$
by using Cauchy's theorem. The possibility of employing this procedure
depends crucially on the analyticity properties of the generating
functions $\mathcal{F}(X,Y)$ (cf. \eqref{def_Fcal}) that, in turn,
hinge upon those of the function $G_{r,r',\sigma}$
(cf.~\eqref{def_ginverse_typeI}).
These can be studied by using the complex implicit-function theorem
\cite{Krantz_02}. Our result can be summarized in the following

\begin{theorem} \label{theo.poly_typeI}
The polynomials  \eqref{def_Pcal} corresponding to 
EGF \eqref{def_Fcal} satisfying Type-I equations are given by
\begin{multline}
\mathcal{P}_n(X)  \;=\;  
\frac{(1+\sigma X)^{n(r-r')+s+s'}}{X^{s+r n}} \\
   \times
   \lim_{Z\rightarrow X}\frac{\partial^n}{\partial Z^n}\left[
   \frac{Z^{s-r-1}}{(1+\sigma Z)^{\eta}}
   \left(\frac{Z-X}{G^{-1}_{r,r',\sigma}(Z)-
                    G^{-1}_{r,r',\sigma}(X)}\right)^{n+1}
   \right]\,,
\label{def_Pcal_typeI}
\end{multline}
where $\eta= s+s'+1+r'-r$, or in the following alternative form if $r\in\Z_0$:
\begin{multline}
\mathcal{P}_n(X) \;=\;  
\frac{(1+\sigma X)^{n(r-r')+s+s'}}{X^{s+r n}}\, \sigma^{(n+1)r}\,
\binom{-1-r'+r}{r}^{-n-1}  \\
\times
   \lim_{Z\rightarrow X}\frac{\partial^n}{\partial Z^n}\left[
   \frac{Z^{s-r-1}(Z-X)^{n+1}}{(1+\sigma Z)^\eta}
   \left[ \log
   \frac{Z\, \widehat{Q}^0_{r,r',\sigma}(Z)}{X\, \widehat{Q}^0_{r,r',\sigma}(X)}
   \right]^{-n-1}
   \right]\,,
\label{def_Pcal_typeI_Bis}
\end{multline}
where
\begin{subequations}
\label{def_Q0hat_Q0}
\begin{align}
\widehat{Q}^0_{r,r',\sigma}(X) &\;=\; \exp \left(
  \frac{Q^{0}_{r,r',\sigma}(X)}
       {\sigma^r \binom{-1-r'+r}{r}}\right) \,, \label{def_Q0hat}\\[2mm] 
       Q^{0}_{r,r',\sigma}(X) &\;=\; \sum\limits_{k\in\Z_0 \setminus\{r\}}
       \sigma^k \, \binom{-1-r'+r}{k} \, \frac{X^{k-r}}{k-r} \,.
\label{def_Q0}
\end{align}
\end{subequations}
\end{theorem}

\proof
Let us pick $X\in\C\setminus\{0\}$ contained within the convergence
disk of $Q^0_{r,r',\sigma}$ ($|X|<1$; cf. \eqref{def_Q0}).
As $Q^0_{r,r',\sigma}$ contains a term of the form $X^{-r}$, the origin can
be a singular point for specific choices of $r$ (either a
pole or a branch point). We consider the function
\be
A \colon U \subset \C^3 \;\rightarrow\; \C \colon (X_1,X_2,X_3)
                        \;\mapsto\; A(X_1,X_2,X_3) \;=\;
   \xi(X_1,X_2)-\xi(X_3,0) \,,
\label{def_function_A}
\ee
where $U \subset \C^3$ is an open neighborhood of $(X,0,X)$ and
\begin{equation}
\xi(X_1,X_2) \;=\; X_2\, X_1^{-r}\, (1+\sigma X_1)^{r-r'}+
G^{-1}_{r,r',\sigma}(X_1)\,.
\label{xi}
\end{equation}
Now, as $A(X,0,X)=0$ and $A_3(X,0,X) \neq 0$ for all $X\in\C$ such that
$0<|X|<1$, there exist open neighborhoods
$U_1\subset \C^2$ and $U_2\subset \C$ of $(X,0)$ and $X$, respectively,
with $U_1\times U_2\subset U$, and a unique holomorphic function
$\theta \colon U_1 \rightarrow U_2$ such that
\be
A^{-1}(0)\cap (U_1\times U_2) \;=\; \left\{ ((X,Y),\theta(X,Y)) \colon
     (X,Y)\in U_1 \right\}\,.
\ee
An important consequence of this result and the definition of $G_{r,r',\sigma}$
(cf. \eqref{def_ginverse_typeI}) is that
\be
G_{r,r',\sigma}(\xi(X,Y)) \;=\; G_{r,r',\sigma}(\xi(\theta(X,Y),0)) \;=\;
\theta(X,Y) \,.
\label{relation.Phi.G}
\ee

The analyticity of $\theta$ implies that there
exists an open neighborhood  $\Omega$ of the origin of the complex $Y$-plane
such that, for every $X\in \mathbb{C}$ satisfying $0<|X|<1$, the function
$Y \mapsto \mathcal{F}(X,Y)$ given by \eqref{FI} is analytic
in $\Omega$. By using Cauchy's theorem we can then write
\be
\hspace*{-2mm}
[Y^n]\mathcal{F}(X,Y) \;=\; \frac{1}{2\pi i}
\int_\Gamma
\left[\frac{G_{r,r',\sigma}(\xi(X,Y))}{X}\right]^s \!
\left[\frac{1+\sigma X}{1+\sigma G_{r,r',\sigma}(\xi(X,Y))}\right]^{s+s^\prime}
\!\!\frac{d Y}{Y^{n+1}} \,,
\label{integral_1}
\ee
where $\Gamma$ is a simple closed curve of index $+1$, contained in $\Omega$,
surrounding the origin $Y=0$ and no other singularity of the
integrand. A natural change of variables, suggested by the form of
\eqref{integral_1}, is to put (cf. \eqref{relation.Phi.G})
\be
Z \;=\; G_{r,r',\sigma}(\xi(X,Y)) \;=\; \theta(X,Y) \,,
\label{def.change_variables}
\ee
so let us consider the one-parameter family of holomorphic maps
($\pi_2$ denotes the projection onto the second argument)
\be
Z_X \colon \pi_2((\{X\}\times \C)\cap U_1)\subset \C \;\rightarrow\;
    \C \colon Y \;\mapsto\; Z_X(Y) \;=\; G(\xi(X,Y))\,.
\ee
Notice that as
$Z_X^\prime(0)\neq0$ we have that $Z_X(\Gamma)$, the image of the
original integration contour $\Gamma$, will be also a closed, simple curve of
index $+1$, contained in $Z_X(\Omega)$ and surrounding the point $Z=X$
in the complex $Z$--plane.
Notice also that, given any open neighborhood
$V_X$ of $Z=X$, it is possible to choose the original integration contour
$\Gamma$ in such a way that $Z_X(\Gamma)\subset V_X$.

It is now straightforward to rewrite the integral in \eqref{integral_1}
as a contour integral on $Z_X(\Gamma)$ to obtain
the following expression for the polynomials $\mathcal{P}_n(X)$:
\begin{multline}
\mathcal{P}_n(X) \;=\;  n! \, \frac{(1+\sigma \, X)^{n(r-r')+s+s'}}
                                 {X^{s+r n}} \\[2mm]
\times \frac{1}{2\pi i} \int_{Z_X(\Gamma)}
\frac{Z^{s-r-1}}{(1+\sigma \, Z)^\eta}
\frac{dZ}{(G^{-1}_{r,r',\sigma}(Z)-G^{-1}_{r,r',\sigma}(X))^{n+1}} \,,
\label{polynomials_calP}
\end{multline}
where $\eta = s+s'+1+r'-r$.

The analytic structure of the integrand of \eqref{polynomials_calP} shows that
it is possible to choose $\Gamma$ in such a way that
(1) the singularity that may appear due to the term $Z^{s-r-1}$ can be avoided,
and (2) the only singularity surrounded by $Z_X(\Gamma)$ is $Z=X$. Hence, we 
can compute the integral by using the residue of the integrand
at $Z=X$. This point can be immediately seen to be a pole of order $n+1$
because $(G^{-1}_{r,r',\sigma})'(X)\neq0$ if $X$ satisfies $0<|X|<1$.
By doing this, we immediately get \eqref{def_Pcal_typeI}.

When $r\in\Z_0$, it is convenient to explicitly take into account the
logarithmic terms appearing in $G^{-1}_{r,r',\sigma}(Z)$ and
$G^{-1}_{r,r',\sigma}(X)$. This gives \eqref{def_Pcal_typeI_Bis}. \hfill \qed

\bigskip

\noindent \textbf{Remarks:}
1. Notice that the procedure that we have followed above
allows us to partially sidestep the difficulties associated with the
impossibility to obtain closed form expressions for the function
$G_{r,r',\sigma}(Z)$ from
\eqref{def_ginverse_typeI} in many cases; we only need the function
$G^{-1}_{r,r',\sigma}(X)$ or $Q^0_{r,r',\sigma}(X)$.

2. As discussed above, the integrand of \eqref{polynomials_calP} has a pole
of order $n+1$ at $Z=X$. The expressions written above are based on the
computation of the residue at this point and take advantage of the fact
that the integrand is a meromorphic function in an open neighborhood of it.
However, in many occasions
it is possible to consider analytic extensions of the integrand
and move the integration contour to rewrite the integral in more convenient
ways.

3. Formula \eqref{def_Pcal_typeI_Bis} suggests the change of variables
$e^U = Z \widehat{Q}^0(Z)$ and $e^V = X \widehat{Q}^0(X)$. This leads to simple
Rodrigues-like formulas for the row polynomials.

%
%
\subsection{Type II case} \label{sec.poly_typeII}

When $\beta\neq 0$ and $\beta'=0$, the EGF $F(x,y)$
has the general form given by Theorem~\ref{theo.typeII}. In this case
we can work with the original variables $x,y$. Our goal is to express
the one-variable polynomials $P_n(x)$ as a contour
integral by following the same steps that led to Theorem~\ref{theo.poly_typeI}.

Using the complex implicit function theorem, we can show (as in the proof
of Theorem~\ref{theo.poly_typeI}) that there is an open neighborhood
$\Omega$ of the origin of the complex $y$-plane where $F(x,y)$
is analytic as a function of $y$ (for every $x\in\C$ satisfying
$0<|x|<1$). Then, $P_n(x)$ can be expressed as a contour integral by using
Cauchy's theorem:
\be
\hspace*{-3mm}
P_n(x)=\frac{n!}{2\pi i} \int_{\Gamma_x}
\left[\frac{G_{\alpha,\beta,\alpha'}(\xi(x,y))}{x}
\right]^{\frac{\alpha+\gamma}{\beta}} \!\!\!\!
\exp\left[\!\!-\frac{\alpha'+\gamma'}{\beta}
\left(x\!-G_{\alpha,\beta,\alpha'}(\xi(x,y))\right)\right]
\frac{dy}{y^{n+1}},
\label{1var_typeII_poly}
\ee
where $\Gamma_x$ is a closed, simple curve of index +1, that surrounds the
origin $y = 0$ and no other singularity of the integrand,
$\xi(x,y)$ is given by
$$
\xi(x,y)\;=\; yx^{-\alpha/\beta}e^{-x\alpha'/\beta}+
              G^{-1}_{\alpha,\beta,\alpha'}(x)\,,
$$
and $G_{\alpha,\beta,\alpha'}$ is defined by \eqref{def_ginverse_typeII}.

As we did in the preceding section, it is convenient now to perform
the change of variables $z=G_{\alpha,\beta,\alpha'}(\xi(x,y))$.
The same steps followed in the proof of Theorem~\ref{theo.poly_typeI} lead to
\be
P_n(x) \;=\;
   \frac{n!}{2\pi i\beta} \,
   \frac{e^{-((n+1)\alpha'+\gamma')x/\beta}}{x^{((n+1)\alpha+\gamma)/\beta}}\,
\int_\Gamma \frac{z^{\gamma/\beta-1} \, e^{\gamma' z/\beta} }
    {\left(G^{-1}_{\alpha,\beta,\alpha'}(z)-
           G^{-1}_{\alpha,\beta,\alpha'}(x)
     \right)^{n+1}}  \, dz\,, \label{typeII_integralformula_general}
\ee
where the integration contour in \eqref{typeII_integralformula_general}
is a closed, simple curve of index +1, that surrounds the point $z=x$
and no other singularity of the integrand. The integrand in
\eqref{typeII_integralformula_general} has a pole of order $n+1$ at $z=x$,
so we can compute this integral using residues. The above discussion can
be summarized in the following

\begin{theorem} \label{theo.poly_typeII}
The polynomials $P_n(x)$ corresponding to EGF satisfying 
Type-II equations are given by
\be
P_n(x) \;=\;
\frac{e^{-((n+1)\alpha'+\gamma')x/\beta}}
     {\beta\, x^{((n+1)\alpha+\gamma)/\beta}}\,
\lim_{z\to x} \frac{\partial^n}{\partial z^n}
\frac{z^{\gamma/\beta-1} \, e^{\gamma' z/\beta}\, (z-x)^{n+1}}
{\left(G^{-1}_{\alpha,\beta,\alpha'}(z)-G^{-1}_{\alpha,\beta,\alpha'}(x)
     \right)^{n+1}} \,.
\label{def_P_typeII}
\ee
\end{theorem}

%
%
\subsection{Type III case} \label{sec.poly_typeIII}

When $\beta=0$ and $\beta'\neq0$, the EGF $F(x,y)$
has the general form given by Theorem~\ref{theo.typeIII}. Again, we
can work with the original variables $x,y$  and write
the polynomials $P_n(x)$ as a contour integral.

Using the complex implicit function theorem, we can show (as in the proof
of Theorem~\ref{theo.poly_typeI}) that there is an open neighborhood
$\Omega$ of the origin of the complex $y$--plane where $F(x,y)$ is analytic
as a function of $y$ (for every $x\in\C$ satisfying $0<|x|<1$). Then,
$P_n(x)$ can be expressed as a contour integral by using Cauchy's theorem:
\be
\hspace*{-3mm}
P_n(x)\!=\!\frac{n!}{2\pi i}\! \int_{\Gamma_x} \!\!
\left[\frac{G_{\alpha,\alpha',\beta'}(\xi(x,y))}{x}
\right]^{1+\frac{\alpha'+\gamma'}{\beta'}}\!\!\!\!\!\!\!
\exp\left[\frac{\alpha+\gamma}{\beta'}
    \left[\frac{1}{x} \!-\!
          \frac{1}{G_{\alpha,\alpha',\beta'}(\xi(x,y))}
\right]\right] \frac{dy}{y^{n+1}},
\label{1var_typeIII_poly}
\ee
where $\Gamma_x$ is a closed, simple curve of index +1, that surrounds the
origin $y = 0$ and no other singularity of the integrand,
$\xi(x,y)$ is given by
$$
\xi(x,y) \;=\; y x^{-\alpha'/\beta'}e^{\alpha/(\beta'x)}+
               G^{-1}_{\alpha,\alpha',\beta'}
$$
and $G_{\alpha,\alpha',\beta'}$ is defined by \eqref{def_ginverse_typeIII}.

As we did in Section~\ref{sec.poly_typeI}, we perform
the change of variables $z=G_{\alpha,\alpha',\beta'}(\xi(x,y))$ to obtain
\be
P_n(x) =
   \frac{n!}{2\pi i \beta'}
   \frac{e^{((n+1)\alpha+\gamma)/(\beta'x)}}
        {x^{1+((n+1)\alpha'+\gamma')/\beta'}}\,
\int _\Gamma \frac{z^{\gamma'/\beta'-1} \, e^{-\gamma/(\beta' z)} }
    {\left(G^{-1}_{\alpha,\alpha',\beta'}(z)-
           G^{-1}_{\alpha,\alpha',\beta'}(x)
     \right)^{n+1}}\, dz\,,
\label{typeIII_integralformula_general}
\ee
where the integration contour in \eqref{typeIII_integralformula_general}
is a closed, simple curve of index +1, that surrounds the point $z=x$
and no other singularity of the integrand. The integrand in
\eqref{typeIII_integralformula_general} has a pole of order $n+1$ at $z=x$.
Then, we can compute this integral using residues.
The above discussion can be summarized in the following

\begin{theorem} \label{theo.poly_typeIII}
The polynomials $P_n(x)$ corresponding to EGF satisfying Type-III equations 
are given by
\be
P_n(x) =
   \frac{e^{((n+1)\alpha+\gamma)/(\beta'x)}}
        {\beta'\, x^{1+((n+1)\alpha'+\gamma')/\beta'}}\,
\lim_{z\to x} \frac{\partial^n}{\partial z^n}
\frac{z^{\gamma'/\beta'-1} \, e^{\gamma/(\beta' z)}\, (z-x)^{n+1}}
    {\left(
     G^{-1}_{\alpha,\alpha',\beta'}(z)-G^{-1}_{\alpha,\alpha',\beta'}(x)
     \right)^{n+1}}
 \,.
\label{def_P_typeIII}
\ee
\end{theorem}

%
%
\subsection{Type IV case} \label{sec.poly_typeIV}

This corresponds to Spivey's case~(S2) \cite{Spivey_11}.
The EGF $F(x,y)$ for the families $(\alpha,0;\alpha',0)$ is given in
closed form by \eqref{def_gf_betas=0}, so it is not necessary to provide
the integral representation used above. The result is easy to
obtain, so we simply quote it here:

\begin{theorem} \label{theo.poly_typeIV}
The polynomials $P_n(x)$ corresponding to EGF satisfying  Type-IV equations 
are given by
\be
P_n(x) \;=\;
   \prod\limits_{k=1}^n \Big(k\, \alpha +\gamma +
                            (k\, \alpha'+ \gamma')\, x \Big)\,.
\label{def_P_typeIV}
\ee
\end{theorem}
Actually, the form of the coefficients for this case is also easy to
obtain:

\begin{corollary} \label{coro.numbers_typeIV}
The coefficients $\binomvert{n}{k}$ for $n\ge 0$ and $0\le k\le n$
corresponding to solutions to Question~\ref{problem.GKP} of Type~IV
are given by
\be
\binomvert{n}{k} \;=\;
   \sum\limits_{t=0}^n \sum\limits_{s=0}^k \stirlingcycle{n}{t} \,
    \binom{t}{s} \, \binom{n-t}{k-s} \,
    (\alpha+\gamma)^{t-s} \, (\alpha'+\gamma')^s
  \alpha^{n-t+s-k} (\alpha')^{k-s} \,.
\label{def_coeff_typeIV}
\ee
\end{corollary}

%
%
\section*{Acknowledgments}

We are indebted to Alan Sokal for his participation in the early stages
of this work, and his encouragement and useful suggestions later on.
We also thank Jesper Jacobsen, Anna de Mier, Neil Sloane, and Mike Spivey
for correspondence, and David Callan for pointing out some interesting
references to us. This research has been supported in part by
Spanish MINECO grant FIS2012-34379. The research of J.S. has also been
supported in part by Spanish MINECO grant MTM2011-24097 and by U.S.\ National
Science Foundation grant PHY--0424082.

\appendix
%
%

\section{Some particular cases}
\label{sec.particular.cases}

In this appendix we will deal with particular cases of interest of the
results obtained in Section~\ref{sec.solutionF}. In particular, we
will give closed formulas for the EGF's corresponding to cases where the
functions $G^{-1}$ of Theorems~\ref{theo.typeI}--\ref{theo.typeIII}
can be written in closed form in terms of simple functions.
Please, note that Type~IV equations have been completely solved in
Theorem~\ref{theo.typeIV}.

%
%
\subsection{Particular cases for Type I equations}
\label{sec.cases.typeI}

In this subsection we will use the
standard parameters $(\alpha,\beta,\gamma;\alpha',\beta',\gamma')$
\eqref{eq_binomvert} in the EGF, instead of the parameters $(r,r',\sigma)$
\eqref{eq_PDE_XY_final}/\eqref{def_rs} labeling  $G^{-1}_{r,r',\sigma}$
(cf.~\eqref{def_ginverse_typeI}). 

%
%
\subsubsection{Solution for $\bm{1+r'=r}$}
\label{sec.spivey3}

This is case (S3) of Spivey \cite{Spivey_11}. The function
$G_{r,r-1,\sigma}^{-1}$ is given by:
\be
G_{r,r-1,\sigma}^{-1}(X) \;=\; \begin{cases}
        \displaystyle - (r \, X^r)^{-1} & \text{if $r\neq 0$,}\\[2mm]
                        \log X          & \text{if $r=0$,}
                \end{cases}
\label{def_Ginverse_spivey3}
\ee
and the EGF is
\be
F(x,y) \;=\; \begin{cases} \displaystyle
\frac{%
\left( \frac{\beta + \beta' \, x \,
       \left( 1 - \alpha \, y (\beta + \beta'\, x)/\beta
       \right)^{-\beta/\alpha}}{\beta + \beta'\, x}
\right)^{\gamma'/\beta' - \gamma/\beta} }{%
    \left( 1 - \frac{\alpha \, y}{\beta} (\beta + \beta' \, x)
             \right)^{(\alpha+\gamma)/\alpha}} & \text{if $\alpha\neq 0$,}
   \\[8mm]
\displaystyle e^{y \, (\beta+\beta' \, x)\gamma/\beta} \,
\left( \frac{\beta + \beta' \, x \, e^{y \, (\beta + \beta' \, x)} }
            {\beta + \beta' \, x}
\right)^{\gamma'/\beta' - \gamma/\beta} & \text{if $\alpha=0$.}
\end{cases}
\label{def_solution_F_final_spivey3}
\ee
This EGF was also obtained by Th\'eor\^et \cite[Proposition~3]{Theoret_95_1}.
A particular case of this family
corresponds to the Eulerian numbers defined by $(0,1,1;1,-1,0)$.

%
%
\subsubsection{Solution for $\bm{r=-1}$}
\label{sec.spivey1}

This is case (S1) of Spivey \cite{Spivey_11}. The function
$G_{-1,r',\sigma}^{-1}$ is given by:
\be
G_{-1,r',\sigma}^{-1}(X) \;=\; \begin{cases}  \displaystyle\frac{\sigma}{(1+r')} \left(
 1 - \frac{1}{(1+\sigma\, X)^{1+r'}} \right)  & \text{if $r'\neq -1$,}\\[3mm]
 \sigma \, \log(1+ \sigma\, X)                & \text{if $r'=-1$,}
                \end{cases}
\label{def_Ginverse_spivey1}
\ee
and this leads to
\be
F(x,y) \;=\; \begin{cases}
\left( \frac{1}{\beta'\, x} \,
  \frac{\beta + \beta' \, x}
       {(1 -(\alpha'+\beta') xy)^{\beta'/(\alpha'+\beta')}}
       - \beta \right)^{\gamma/\beta -1} & \text{} \\[4mm]
\qquad  \qquad \times
\left( \frac{1}{ 1 -(\alpha'+\beta') xy}
\right)^{ \frac{2\beta\, \beta' + (\alpha'+\gamma')\, \beta -\gamma\, \beta'}
  {\beta\, (\alpha'+\beta')} } & \text{if $\alpha'\neq -\beta'$,} \\[5mm]
e^{x\, y\, \beta'\, (1+\gamma'/\beta' - \gamma/\beta)} \,
\left( \frac{ e^{x\, y\, \beta'} (\beta + \beta' \, x) - \beta}{\beta' \, x}
\right)^{-1+\gamma/\beta} & \text{if $\alpha'=-\beta'$.}
\end{cases}
\label{def_solution_F_final_spivey1}
\ee
This EGF can also be obtained starting from the EGF for $r'=0$
\eqref{def_solution_F_final_rp=0}, and using the involution
\eqref{def_duality} (cf. first row of Table~\ref{table.involutions}),
as suggested by Th\'eor\^et \cite[p.~97]{Theoret_95_1}.

%
%
\subsubsection{Solution for $\bm{r=r'=1}$}
\label{sec.r=rp=1}

In this case, the function to invert is
\be
G_{1,1,\sigma}^{-1}(X) \;=\;  -1/X +
        \sigma \, \log \left( (1+\sigma\, X)/X \right)\,.
\ee
This immediately gives
\begin{multline}
F(x,y) \;=\;  
\frac{1}{(\beta' x)^{1+\gamma/\beta}} \,
\left( \frac{
    T\left( e^{\beta^2\, y/(\beta'\, x)} \, T^{-1}(1+\beta/(\beta' x))\right)}
    {\beta+\beta'\, x} \right)^{1+\gamma'/\beta'-\gamma/\beta} \\
\times \left( \frac{ \beta }
{T \left( e^{\beta^2\, y/(\beta' \, x)} \, T^{-1}(1+\beta/(\beta' x))
\right) -1} \right)^{2+\gamma'/\beta'} \,,
\label{def_solution_F_final_r=rp=1}
\end{multline}
in terms of the tree function $T$ \cite{Corless_96,Corless_97}.
An interesting particular case corresponds to $(1,1,-1;1,1,-1)$
giving the numbers $\associatedstirlingcycleBis{n+k}{k}$ introduced
in \cite{Fekete_94}.

%
%
\subsubsection{Solution for $\bm{r'=0}$}
\label{sec.rp=0}

This is the case studied by Neuwirth \cite{Neuwirth_01}.
The function $G_{r,0,\sigma}^{-1}$ is given by:
\be
G^{-1}_{r,0,\sigma}(X) \;=\; \begin{cases}
   \displaystyle -(1/r)\, \left( (1+\sigma\, X)/X\right)^r
   & \text{if $r\neq 0$,}\\[2mm]
 \displaystyle  \log \left( X/(1 + \sigma\, X)\right)  & \text{if $r=0$.}
                \end{cases}
\label{def_Ginverse_rp=0}
\ee
Then, we have
\be
F(x,y) \;=\; \begin{cases}\displaystyle
\left( \frac{\beta}
   {(\beta+\beta'\, x)(1-\alpha\, y)^{\beta/\alpha} -\beta'\, x}
             \right)^{1+ \gamma'/\beta'} &  \text{}  \\[4mm]
\qquad\qquad \times
  (1- \alpha \, y)^{ (\beta/\alpha)(1+\gamma'/\beta'-(\alpha+\gamma)/\beta)} &
\text{if $\alpha\neq 0$,} \\[5mm]
\displaystyle
e^{y \, \gamma}\, \left(
   \frac{ \beta }{\beta + \beta'\, x \, (1 - e^{y \, \beta})}
   \right)^{1 + \gamma'/\beta'} & \text{if $\alpha=0$.}
\end{cases}
\label{def_solution_F_final_rp=0}
\ee
This EGF was also obtained by Th\'eor\^et
\cite{Theoret_thesis,Theoret_95_1}.\footnote{
  The expression for the case $\alpha=0$ has a typo in
  \cite[Eq.~(4.66)]{Theoret_thesis}; but it is correct in
  \cite[Eq.~(16)]{Theoret_95_1}.
}
An important particular case corresponds to the numbers
$\mathrm{Surj}(n,k)$ \cite{Fekete_94} defined by $(0,1,0;0,1,0)$.

%
%
\subsubsection{Solution for $\bm{r=0}$ and $\bm{-r'\in \N}$}
\label{sec.r=0_rp=-nu}

As $r'$ is a negative integer, it is convenient to define $\nu = -r' \in\N$.
Then, the function $G_{0,-\nu,\sigma}^{-1}$ is given by
\be
G_{0,-\nu,\sigma}^{-1}(X) \;=\; \log X + \sum\limits_{k=1}^{\nu-1}
\binom{\nu-1}{k} \, \frac{(\sigma\, X)^k}{k} \,.
\ee
It is convenient to define a new function $T_\nu$ as
\be
T_\nu^{-1}(z) \;=\; z \, e^{ Q_\nu(z) } \,, \quad \text{where} \quad
Q_\nu(z)      \;=\; \sum\limits_{k=1}^{\nu-1} \binom{\nu-1}{k} \,
                    \frac{(-z)^k}{k} \,.
\label{def_Tnu}
\ee
It is clear that $T_1$ is the identity, and $T_2$ is the tree function
\cite{Corless_96,Corless_97}. The EGF is
\begin{multline}
F(x,y) \;=\;  \beta^{1-\nu+\gamma'/\beta'}\,
\left(
\frac{ T_\nu \left( e^{y\, \beta^{1-\nu}\, (\beta+\beta'\,x)^\nu} \,
       T_\nu^{-1}(-\beta'\, x/\beta) \right) }{ (-\beta' \, x)}
\right)^{\gamma/\beta} \\
 \times \; \left(
\frac{1- T_\nu \left( e^{y\, \beta^{1-\nu}\, (\beta+\beta'\,x)^\nu} \,
       T_\nu^{-1}(-\beta'\, x/\beta) \right) }{\beta + \beta'\, x}
\right)^{1-\nu-\gamma/\beta+\gamma'/\beta'} \,.
\label{def_solution_F_final_r=0_rp=-nu}
\end{multline}
It is possible to define $\nu$-order Eulerian numbers as a
generalization of ordinary and second order Eulerian numbers
by the parameter choice $(0,1,1;\nu,-1,1-\nu)$ \cite{BSV}.

%
%
\subsubsection{Solution for $\bm{r=0}$ and $\bm{r'\in \N}$}
\label{sec.r=0_rp=nu}

It is convenient to redefine $r'=\nu\in\N$ in accordance with the previous
section. The function $G_{0,\nu,\sigma}^{-1}$ is given by
\begin{equation}
G_{0,\nu,\sigma}^{-1}(X)=\log\left( \sigma \, T_{\nu+1}^{-1} \left(
                  \frac{\sigma\, X}{1+\sigma\, X} \right) \right)\,,
\end{equation}
where $T_\nu$ is given by \eqref{def_Tnu}. We have now that
\begin{multline}
F(x,y) \;=\;  \frac{\beta^{1+\nu+\gamma'/\beta'}}
       {(\beta + \beta'\, x)^{1+\nu+\gamma'/\beta'-\gamma/\beta}}\,
\left(
\frac{ T_{\nu+1} \left( e^{y\, \beta^{1+\nu}\, (\beta+\beta'\,x)^{-\nu}} \,
       T_{\nu+1}^{-1}(\beta'\, x/(\beta+\beta'\,x)) \right) }{ \beta' \, x}
\right)^{\gamma/\beta} \\
\times \left(
1- T_{\nu+1} \left( e^{y\, \beta^{1+\nu}\, (\beta+\beta'\,x)^{-\nu}} \,
       T_{\nu+1}^{-1}(\beta'\, x/(\beta+\beta'\,x)) \right)
\right)^{-1-\nu-\gamma'/\beta'} \,.
\label{def_solution_F_final_r=0_rp=nu}
\end{multline}
Interesting particular cases are the $\nu$-order Ward numbers
--a generalization of the ordinary Ward numbers--
given by $(0,1,0;\nu,1,-\nu)$ \cite{BSV}.

%
%
\subsection{Particular cases for Type II equations}
\label{sec.cases.typeII}

In this subsection we will use the standard parameters
$(\alpha,\beta,\gamma;\alpha',0,\gamma')$ in both
$G_{\alpha,\beta,\alpha'}^{-1}$ (cf.~\eqref{def_ginverse_typeII}) and the EGF.

%
%
\subsubsection{Solution for $\bm{(-\beta,\beta;\alpha',0)}$}
\label{sec.a=-b_bp=0}

In this case $G_{-\beta,\beta,\alpha'}^{-1}$ can be
computed in closed form to give
\be
G_{-\beta,\beta,\alpha'}^{-1}(x) \;=\; \begin{cases}
                \left(1-e^{-\alpha' \, x/\beta}\right)/\alpha' &
                            \text{if $\alpha'\neq 0$,}  \\[2mm]
                 x/\beta  & \text{if $\alpha'= 0$.}
 \end{cases}
\ee
Hence
\be
F(x,y) \;=\;
\begin{cases}
   \displaystyle
   \left(1 - \frac{\beta \, \log(1 - \alpha' \, x\, y)}{\alpha'\,x}
             \right)^{\gamma/\beta-1} \,
            (1 - \alpha' \, x\, y)^{-1-\gamma'/\alpha'}
    &\text{if $\alpha'\neq0$,} \\[5mm]
   \displaystyle
    (1+\beta \, y)^{-1+\gamma/\beta} e^{\gamma' \, x\, y } &
     \text{if $\alpha'=0$.}
\label{def_gf_a=-b_bp=0}
\end{cases}
\ee
A relevant particular case is defined by the parameters $(1,-1,0;0,0,1)$ and
corresponds to the numbers $n^{\underline{n-k}}$.

%
%
\subsubsection{Solution for $\bm{(\alpha,\beta;0,0)}$}
\label{sec.ap=0_bp=0}

In this case $G_{\alpha,\beta,0}^{-1}$ can be easily summed to give
\be
G_{\alpha,\beta,0}^{-1}(x) \;=\; \begin{cases}
     - x^{-\alpha/\beta}/\alpha & \text{if $\alpha\neq 0$,} \\[2mm]
     (1/\beta) \, \log x        & \text{if $\alpha = 0$.}
 \end{cases}
\ee
Hence we have
\be
F(x,y) \;=\; \begin{cases}
   (1- \alpha \, y)^{-(1 + \gamma/\alpha)} \, \exp\left(
    -\frac{\gamma'\, x}{\beta}\, (1 - (1-\alpha\, y)^{-\beta/\alpha} ) \right)
    &\text{if $\alpha\neq0$,} \\[2mm]
    \exp\left(\gamma \, y-\gamma'\, (1-e^{\beta \, y})x/\beta\right)
    &\text{if $\alpha=0$.}
  \end{cases}
\label{def_gf_betap=0_alphap=0}
\ee
The generalization of the Lah and Stirling subset numbers
$S(r;n,k)$ \cite{Lang_00,Lang_09} corresponding to $(r-1,1,1-r;0,0,1)$ is 
a particular interesting case.

%
%
\subsection{Particular cases for Type III equations}
\label{sec.cases.typeIII}

In this subsection we will use the standard parameters
$(\alpha,0,\gamma;\alpha',\beta',\gamma')$ in both
$G_{\alpha,\alpha',\beta'}^{-1}$ (cf.~\eqref{def_ginverse_typeIII}) and the EGF.

%
%
\subsubsection{Solution for $\bm{(\alpha,0;\beta',\beta')}$}
\label{sec.ap=bp_b=0}

The expression for $G_{\alpha,\beta',\beta'}^{-1}$ is given by
\be
G_{\alpha,\beta',\beta'}^{-1}(x) \;=\; \begin{cases}
    \displaystyle
    e^{\alpha/(\beta'\, x)}\, \left( \beta'/\alpha^2 - 1/(\alpha\, x) \right)
      -\beta'/\alpha^2 & \text{if $\alpha\neq 0$,} \\[2mm]
      -1/(2\beta' \, x^2) & \text{if $\alpha=0$.}
\end{cases}
\ee
Therefore,
\be
F(x,y) = \begin{cases}
 \left[\frac{\alpha}{x\, \beta' \, (1- T(\zeta))} \right]^{2+\gamma'/\beta'}\,
 \exp\left[ \frac{\alpha+\gamma}{\beta'}\,
 \left(\frac{1}{x} - \frac{\beta'}{\alpha}\, (1-T(\zeta)) \right)\right]
 & \text{if $\alpha\neq0$,} \\[5mm]
 (1-2\beta' x y)^{-(1 + \gamma'/(2\beta'))} \,
 \exp\left[ \frac{\gamma}{\beta' x}\, ( 1 - \sqrt{1- 2 \beta' x y} )
 \right]
 & \text{if $\alpha=0$,}
\end{cases}
\label{def_gf_beta=0_alphap=betap}
\ee
where
\be
\zeta(x,y) \;=\;
 \left( \alpha^2 \, y + \beta' \, x - \alpha\right) \,
     \exp\left(\alpha/(\beta' \, x)-1\right)/(\beta' \, x)  \,.
\ee
The particular family $(1,0,\gamma;1,1,\gamma')$ contains numbers
related to the Ramanujan functions $Q_{n,k}(x)$ \cite{Zeng_99} (see
Table~\ref{table.example}). These EGF's seem to be new.

%
%
\subsubsection{Solution for $\bm{(0,0;\alpha',\beta')}$}
\label{sec.a=0_b=0}

A straightforward computation leads to the expression for
$G_{0,\alpha',\beta'}^{-1}$:
\be
G_{0,\alpha',\beta'}^{-1}(x) \;=\; \begin{cases}
 -x^{-1-\alpha'/\beta'}/(\alpha'+\beta')
    & \text{if $\alpha'\neq -\beta'$,} \\[2mm]
  (\log x)/\beta'  & \text{if $\alpha'=-\beta'$.}
\end{cases}
\ee
Hence
\be
F(x,y) \;=\; \begin{cases}
 \displaystyle \frac{ \exp \Big( \gamma\,
  \big(1 - (1-x \, y\, (\alpha'+\beta'))^{\beta'/(\alpha'+\beta')}
  \big)/(\beta'\, x)
  \Big)}
  {\big(1- x \, y \, (\alpha'+\beta')\big)^{1 + \gamma'/(\alpha'+\beta')}}
 & \text{if $\alpha'\neq-\beta'$,} \\[5mm]
e^{x\, y\, \gamma'}
 \, \exp\left( \gamma\,
  \left(1 - e^{-x\, y\, \beta'} \right) /(\beta'\, x)\right)
 & \text{if $\alpha'=-\beta'$.}
\end{cases}
\label{def_gf_alpha=0_beta=0}
\ee

%
%
\subsubsection{Solution for $\bm{(\alpha,0;0,\beta')}$}
\label{sec.ap=0_b=0}

The expression for $G_{\alpha,0,\beta'}^{-1}$ is given by:
\be
G_{\alpha,0,\beta'}^{-1}(x) \;=\; \begin{cases}\displaystyle
    (1-e^{\alpha/(\beta' x)})/\alpha  & \text{if $\alpha\neq 0$,} \\[2mm]
    - 1/(\beta'\, x) & \text{if $\alpha=0$.}
\end{cases}
\ee
The corresponding EGF is
\be
F(x,y) \;=\; \begin{cases}\displaystyle
\left( \frac{\alpha}{\alpha + \beta' \, x\, \log(1-\alpha\, y)}
             \right)^{1+\gamma'/\beta'} \,
 (1- \alpha \, y)^{-(1+\gamma/\alpha)}
 & \text{if $\alpha\neq0$,} \\[5mm]
\displaystyle e^{\gamma\, y} \, (1- x\, y\, \beta')^{-(1+\gamma'/\beta')}
 & \text{if $\alpha=0$.}
\end{cases}
\label{def_gf_alphap=0_beta=0}
\ee
A particular case corresponds to the injective numbers
$\mathrm{Inj}(n,k)$ \cite{Fekete_94}
defined by the parameters $(0,0,1;0,1,0)$.

%
%

\end{document}